\documentclass[a4paper,12pt,reqno]{amsart}
\usepackage{graphics}
\newcommand{\mj}{\operatorname{MJ}}
\newcommand{\fj}{\operatorname{FJ}}
\newcommand{\bj}{\operatorname{BJ}}

\newcommand{\id}{\operatorname{id}}

\newcommand{\switch}{\operatorname{SWITCH}}

\newtheorem{defi}{Definition}
\newtheorem{theo}{Theorem}

\newtheorem{cor}{Corollary}

\newtheorem{ex}{Example}

\begin{document}

\author[Ilse Fischer]{\box\Adr}

\newbox\Adr
\setbox\Adr\vbox{
\centerline{ \large Ilse Fischer}
\vspace{0.3cm}
\centerline{Institut f\"ur Mathematik der Universit\"at Klagenfurt,}
\centerline{Universit\"atsstrasse 65-67, A-9020 Klagenfurt, Austria.}
\centerline{E-mail: {\tt Ilse.Fischer@uni-klu.ac.at}}
}

\title[Symmetry theorem]{A symmetry theorem on a modified jeu de taquin}

\date{}

\begin{abstract}
For their bijective proof of the hook-length formula for the number of
standard tableaux of a fixed shape Novelli, Pak and Stoyanovskii 
\cite{NPS} define a modified jeu de taquin which transforms an arbitrary
filling of the Ferrers diagram with $1,2,\dots,n$ (tabloid) into a standard
tableau. Their definition relies on a total order of the cells in the Ferrers
diagram induced by a special standard tableau, however, this definition also
makes sense for the total order induced by any other standard tableau. Given two
standard tableaux $P,Q$ of the same shape we show that the number of tabloids
which result in $P$ if we perform modified jeu de taquin with respect to the
total order induced by $Q$ is equal to the number of tabloids which result in
$Q$ if we perform modified jeu de taquin with respect to $P$. This symmetry
theorem extends to skew shapes and shifted skew shapes.
\end{abstract}   

\maketitle

\section{Introduction}

\begin{figure}
\setlength{\unitlength}{1cm}
\begin{picture}(14.5,4)
\put(1,1){\scalebox{0.35}{\includegraphics{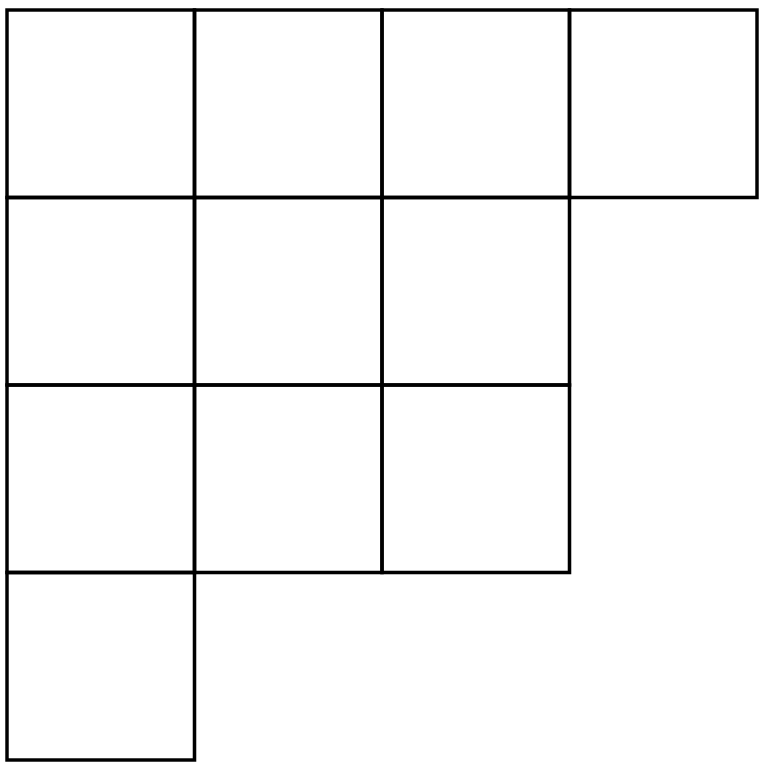}}}
\put(8,1){\scalebox{0.35}{\includegraphics{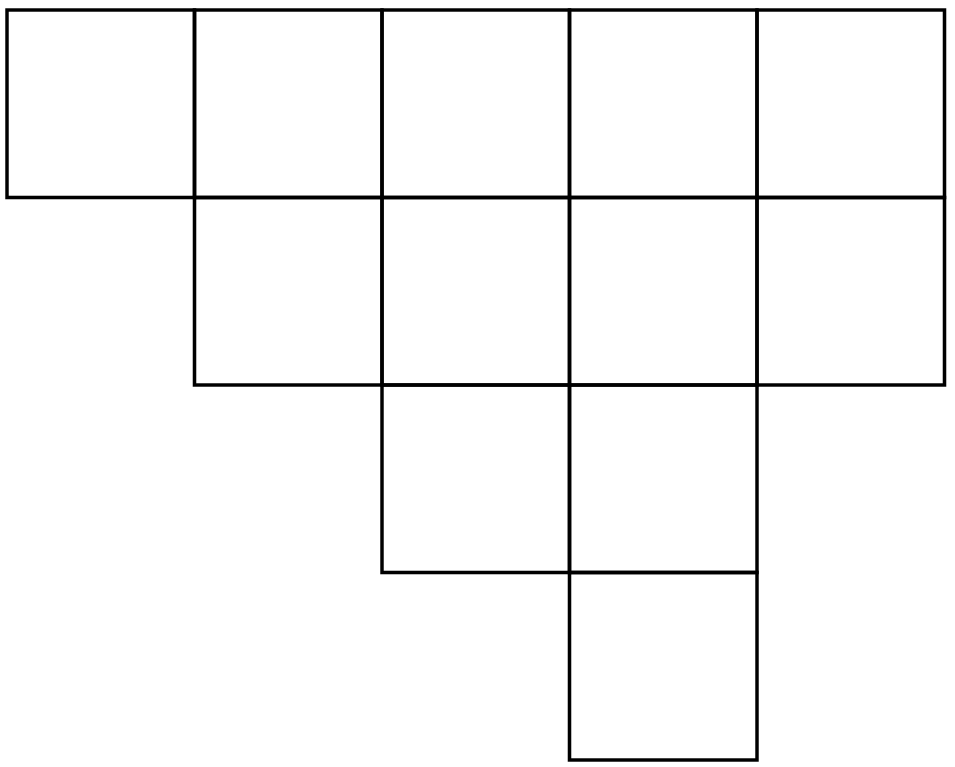}}}
\put(1,0.5){a. The Ferrers diagram}
\put(1,0){corresponding to $(4,3,3,1)$}
\put(8,0.5){b. The shifted Ferrers diagram}
\put(8,0){corresponding to $(5,4,2,1)$} 
\end{picture}
\caption{}
\label{ferrer}   
\end{figure}

\begin{figure}
\setlength{\unitlength}{1cm}
\begin{picture}(14.5,4)
\put(1,1){\scalebox{0.35}{\includegraphics{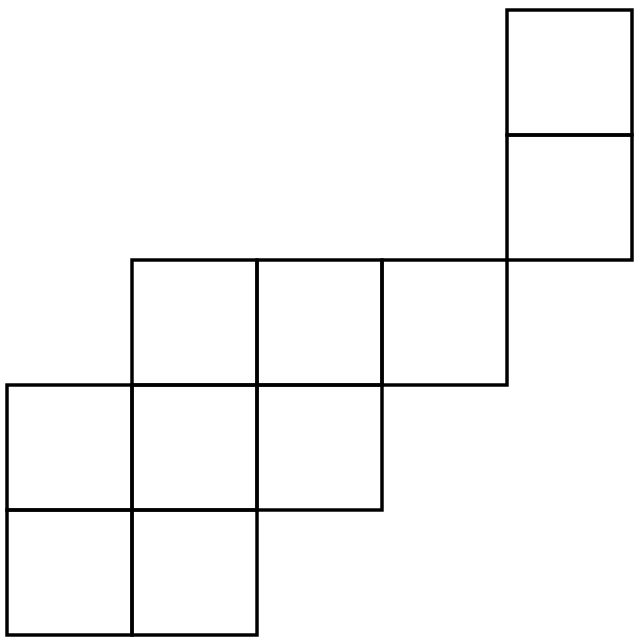}}}
\put(8,1){\scalebox{0.35}{\includegraphics{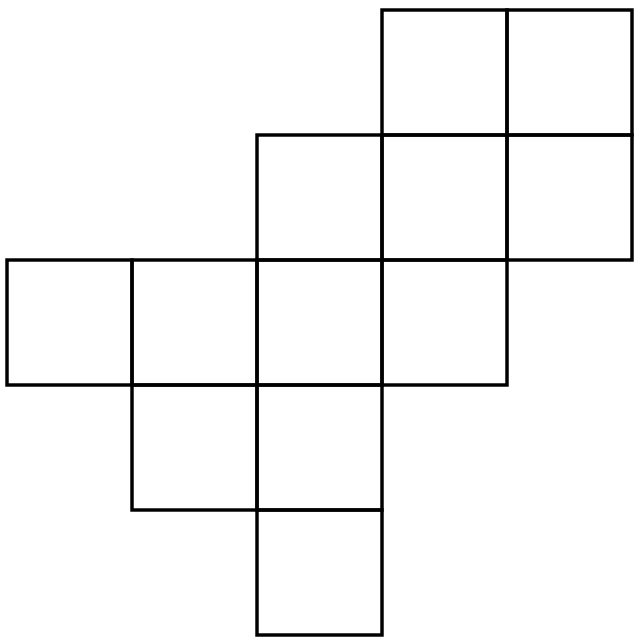}}}
\put(1,0.5){a. The skew Ferrers diagram of}
\put(1,0){shape $(5,5,4,3,2) / (4,4,1)$.}
\put(8,0.5){b. The shifted skew Ferrers diagram of}
\put(8,0){shape $(7,6,4,2,1) / (5,3)$.} 
\end{picture}
\caption{}
\label{skew}   
\end{figure}

A partition of a positive integer $n$ is  a sequence of integers
$\lambda=(\lambda_{1},\lambda_{2},\dots,\lambda_{r})$ with
$\lambda_{1}+\lambda_{2}+\dots+\lambda_{r}=n$
and $\lambda_{1}\ge \lambda_{2}\ge \dots \ge \lambda_{r} \ge 0$. The {\it (unshifted) Ferrers
diagram} of shape $\lambda$ is an array of cells with $r$ left-justified rows 
and $\lambda_{i}$ cells in row $i$. (See Figure~\ref{ferrer}.a.) If $\lambda$
is a partition with distinct components 
then the {\it shifted Ferrers diagram} of shape $\lambda$ is an array of
cells with $r$ rows, each row indented by one cell to the right with respect
to the previous row and $\lambda_{i}$ cells in row $i$. (See
Figure~\ref{ferrer}.b.) If 
$\mu = (\mu_{1},\mu_{2}, \dots, \mu_{r})$,
$\lambda=(\lambda_{1},\lambda_{2},\dots,\lambda_{r})$ 
are partitions (resp. partitions with distinct components) such
that $\lambda_{i} \le \mu_{i}$ for $1 \le i \le r$ then the unshifted (resp. shifted) skew 
Ferrers diagram of shape $\mu / \lambda$ is the diagram we obtain if we remove the cells of the
unshifted (resp. shifted) Ferrers diagram of shape $\lambda$ from the
unshifted (resp. shifted) Ferrers diagram of shape $\mu$. 
(See Figure~\ref{skew}.)

\medskip

\begin{defi}  
A tabloid of (shifted) skew shape $\mu / \lambda$ is an arbitrary
filling of the (shifted) skew Ferrers diagram of shape 
$\mu / \lambda$ with the integers $1, 2, \dots, |\mu / \lambda|$. 
\end{defi}

\medskip

As usual a (shifted) standard skew tableau is a tabloid with increasing rows and columns.
For the rest of the article we fix an unshifted or shifted skew shape $\mu / \lambda$ and set 
$n = | \mu / \lambda|$.

\medskip

\begin{ex}
\label{ex1}
Observe that  
$$R= \begin{array}{cccccc} &   &   &   &   & 8  \\  
                           &   &   &   & 3 & 6  \\ 
                           &   & 9 & 5 & 1 & 4  \\
                           &   &   & 2 & 7 &      
\end{array}$$
is a tabloid of shifted skew shape $(6,5,4,2) / (5,3)$ and that 
$$P= \begin{array}{cccccc} &   &   &   &   & 2  \\  
                           &   &   &   & 1 & 5  \\ 
                           &   & 3 & 4 & 6 & 8  \\
                           &   &   & 7 & 9 &       
\end{array}$$
is a shifted standard skew tableau of shifted skew shape $(6,5,4,2) / (5,3)$.
\end{ex}

\medskip

Next we define forward jeu de taquin in a tabloid $T$ of (shifted) skew shape $\mu /
\lambda$. 
If $(i,j)$ is a cell in the (shifted) skew Ferrers diagram of shape $\mu / \lambda$ let $T_{i,j}$
denote its entry in $T$. In order to simplify the description
of forward jeu de taquin we define $T_{i,j}=\infty$ if
$(i,j)$ is a cell outside of $T$.  

\begin{defi} [Forward jeu de taquin]
Let $T$ be a tabloid and $e$ an entry in $T$. Forward jeu de
  taquin in $T$ with $e$ is defined as follows: Consider the neighbour of $e$
  to the right and the neighbour of $e$ below and 
  if $e$ is greater than the minimum of these neighbours we exchange $e$
  with the minimum. Next consider the new neighbours of $e$ to the right and below and
exchange $e$ with the minimum of these two if $e$ is greater than this
minimum. We repeat this procedure with $e$ until $e$ is stable, i.e. smaller
than its neighbour to the right and its neighbour below. 
\end{defi}

\smallskip

\begin{ex}
\label{fjex}
Performing forward jeu de taquin with $8$ in the tabloid 
$$\begin{array}{ccc} 8 & 1 & 4 \\ 2 & 3 & 5 \\ 6 & 7 & \end{array}$$
results in the sequence 
$$\begin{array}{ccc} 1 & 8 & 4 \\ 2 & 3 & 5 \\ 6 & 7 & \end{array}, \qquad 
  \begin{array}{ccc} 1 & 3 & 4 \\ 2 & 8 & 5 \\ 6 & 7 & \end{array}, \qquad 
  \begin{array}{ccc} 1 & 3 & 4 \\ 2 & 5 & 8 \\ 6 & 7 & \end{array}. \qquad $$
\end{ex}  

\medskip

Next we define an 'ordering procedure' which assigns a
(shifted) standard skew tableau to every tabloid $T$. The procedure depends on another (shifted)
standard skew tableau $S$ of the same shape.

\begin{defi}[Modified jeu de taquin]
\label{mjt}
Let $T$ be a tabloid and $S$ be a (shifted) standard skew tableau of the same shape. Modified 
jeu de taquin 
in $T$ with respect to $S$ is defined as the step by step performance of forward jeu de taquin 
with the entries in 
$T$ in the order the (shifted) standard skew tableau $S$ predicts, starting 
with the entry in $T$, whose cell has the greatest label in $S$.  We denote the resulting 
(shifted) standard skew tableau by $\mj_S(T)$.
\end{defi}

\smallskip

Observe that for tabloids of normal shape $(3,3,2)$ the 'ordering procedure' 
from \cite{NPS} is the modified jeu de taquin with respect to the 
standard tableau 
$\begin{array}{ccc} 1 & 4 & 7 \\ 2 & 5 & 8 \\ 3 & 6 & \end{array}$ from
Definition~\ref{mjt}. In their paper Novelli, Pak and Stoyanovskii show that their 
'ordering procedure' has the nice property that the number of
tabloids that are mapped to a fixed standard tableau is independent of this 
output standard tableau and that this number is equal to the product over all hook-lengths. 

\medskip

In our running example (Example~\ref{ex1}) modified jeu de taquin with respect
to $P$ applied to the tabloid $R$ gives the intermediate tabloids 
$$   
\begin{array}{cccccc} &   &   &   &   & 8  \\  
                           &   &   &   & 3 & 4  \\ 
                           &   & 9 & 5 & 1 & 6  \\
                           &   &   & 2 & 7        
\end{array}, \quad 
     \begin{array}{cccccc} &   &   &   &   & 8  \\  
                           &   &   &   & 3 & 4  \\ 
                           &   & 9 & 1 & 5 & 6  \\
                           &   &   & 2 & 7        
\end{array}, \quad 
     \begin{array}{cccccc} &   &   &   &   & 8  \\  
                           &   &   &   & 3 & 4  \\ 
                           &   & 1 & 2 & 5 & 6  \\
                           &   &   & 7 & 9        
\end{array},
$$
before we finally obtain
$$\mj_P(R)= 
     \begin{array}{cccccc} &   &   &   &   & 4  \\  
                           &   &   &   & 3 & 6  \\ 
                           &   & 1 & 2 & 5 & 8  \\
                           &   &   & 7 & 9        
\end{array}
 =:Q.$$

\section{The symmerty theorem}
\label{symsec}

If $P,Q$ are two (shifted) standard skew tableaux of the same shape $\mu / \lambda$ let $A_{P, Q}$ denote 
the number of tabloids $T$ of shape $\mu / \lambda$ with the property that the application of modified jeu 
de taquin to 
$T$ with respect to $P$ results in  $Q$.  For the normal shape
$(3,3,2)$ the matrix $(A_{P,Q})_{(P,Q)}$ is presented in Figure~\ref{matrix},
where  '$1$' stands for $936$, '$2$' stands for $944$, 
'$3$' stands for $960$, '$4$' stands for $976$, 
'$5$' stands for $984$ and '$6$' stands for $996$. Moreover 
every row and column corresponds to one of the $42$ standard 
tableaux of shape $(3,3,2)$, which are ordered lexicographically if we
identify a standard tableau with the permutation we obtain by reading the standard
tableau rowwise from top to bottom and within a row from left to right.

\begin{figure}
{\tiny
\begin{eqnarray*}
333333333333333333333333333333333333333333\\
333333333333333333333333333333333333333333\\
333333333333333333333333333333333333333333\\
333333333333333333333333333333333333333333\\
333333333333333333333333333333333333333333\\
333334445522255222222333333333333333333333\\
333334445522255222222333333333333333333333\\
333334445522255222222333333333333333333333\\
333335556611166111111333333333333333333333\\
333335556611166111111333333333333333333333\\
333332221144411444444333333333333333333333\\
333332221144411444444333333333333333333333\\
333332221144411444444333333333333333333333\\
333335556611166111111333333333333333333333\\
333335556611166111111333333333333333333333\\
333332221144411444444333333333333333333333\\
333332221144411444444333333333333333333333\\
333332221144411444444333333333333333333333\\
333332221144411444444333333333333333333333\\
333332221144411444444333333333333333333333\\
333332221144411444444333333333333333333333\\
333333333333333333333444114441144224233333\\
333333333333333333333444114441144224233333\\
333333333333333333333444114441144224233333\\
333333333333333333333111661116611551533333\\
333333333333333333333111661116611551533333\\
333333333333333333333444114441144224233333\\
333333333333333333333444114441144224233333\\
333333333333333333333444114441144224233333\\
333333333333333333333111661116611551533333\\
333333333333333333333111661116611551533333\\
333333333333333333333444114441144224233333\\
333333333333333333333444114441144224233333\\
333333333333333333333222552225522442433333\\
333333333333333333333222552225522442433333\\
333333333333333333333444114441144224233333\\
333333333333333333333222552225522442433333\\
333333333333333333333333333333333333333333\\
333333333333333333333333333333333333333333\\
333333333333333333333333333333333333333333\\
333333333333333333333333333333333333333333\\
333333333333333333333333333333333333333333\\
\end{eqnarray*}}
\caption{The Matrix $(A_{P,Q})_{(P,Q)}$ for the normal shape $(3,3,2)$.}
\label{matrix}
\end{figure}

By this and other computer experiments\footnote{Those computer experiments
  were originally intended to find a total order of the
cells in the shifted Ferrers diagram such that the number of tabloids that are
mapped to a fixed shifted standard tableau by modified jeu de taquin with
respect to the order is independent of the output
shifted standard tableau and with this a proof of the shifted hook-length
formula similar to \cite{NPS}. In \cite{shook} we show that the rowwise performance
of modified jeu de taquin, from bottom to top and within a row from right to
left, has this nice property in the shifted case.} we were led to the conjecture
that $A_{P,Q}=A_{Q,P}$. However, we discovered that a more general theorem is the 
key to this observation. In order to state it, we need another definition.

\medskip

\begin{defi}
Let $S$ be a (shifted) standard skew tableau and $\pi$ a permutation of $\{1,2,\dots,n\}$. Then $S_\pi$
denotes the tabloid we obtain from $S$ by replacing every entry $i$ in $S$ 
by $\pi(i)$.
\end{defi}

\medskip

Let $\pi = 3 \, 8 \, 9 \, 5 \, 6 \, 1 \, 2 \, 4 \, 7$. Observe that $R =P_\pi$
and
$$
Q_{\pi^{-1}} =    \begin{array}{cccccc} &   &   &   &   & 8  \\  
                           &   &   &   & 1 & 5  \\ 
                           &   & 6 & 7 & 4 & 2  \\
                           &   &   & 9 & 3        
\end{array}.$$  
The fact that $\mj_Q(Q_{\pi^{-1}}) = P$ does not come by chance.

\medskip

\begin{theo} 
\label{sym}
Let $P, Q$ be two (shifted) standard skew tableaux of the same 
shape and let $\pi$ be a permutation of $\{1,2,\dots,n\}$. Then
$$\mj_P(P_\pi) = Q \Leftrightarrow \mj_Q(Q_{\pi^{-1}}) = P.$$
\end{theo}

\smallskip

Before we are able to prove the theorem we need the definition of backward jeu
de taquin which is in some sense the inverse operation of forward jeu de taquin. In order to
simplify the description we set $T_{i,j}=0$ if $(i,j)$
is a cell outside of $T$, where $T$ is a tabloid.

\begin{defi} [Backward jeu de taquin]
\label{bj}
Let $T$ be a tabloid and $e$ an entry in $T$. Backward jeu de taquin in $T$
with $e$ is defined as follows: Exchange $e$ with the maximum of its neighbour
to the left and its neighbour above if $e$ has either a neighbour to the left
or above in the fixed shape. We repeat this procedure with $e$ until
$e$ has no neighbour to the left and no neighbour above in the fixed (shifted)
skew shape.
\end{defi}

\smallskip

Observe that in Example~$\ref{fjex}$ the input tabloid can be obtained from
the output tabloid by performing backward jeu de taquin with $8$.
 
\medskip

{\it Proof of Theorem~\ref{sym}.} We only have to show one direction of the assertion, for the other follows by 
symmetry. 

Let $S$ be a (shifted) standard skew tableau and $T$ be a tabloid. We define
$\fj(T,S)=(S',T')$, where $T'=\mj_S(T)$ and $S'$ is the tabloid 
we obtain from $S$ after the performance of modified jeu de taquin in $T$ with respect to $S$, 
if we simultaneously apply the transpositions we apply during modified jeu de taquin  
to $T$ also to $S$. (If we exchange $T_{i,j}$ and $T_{i,j+1}$ in $T$ we
exchange $S_{i,j}$ and $S_{i,j+1}$ in $S$ and if we exchange $T_{i,j}$ and
$T_{i+1,j}$ in $T$ we exchange $S_{i,j}$ and $S_{i+1,j}$ in $S$.)

In our running example: If we perform modified jeu de taquin in $R$ with respect to
$P$ and perform the transpositions simultaneously in $P$ we obtain 
the intermediate tabloids
$$
     \begin{array}{cccccc} &   &   &   &   & 2  \\  
                           &   &   &   & 1 & 8  \\ 
                           &   & 3 & 4 & 6 & 5  \\
                           &   &   & 7 & 9        
\end{array}, \quad
\begin{array}{cccccc} &   &   &   &   & 2  \\  
                           &   &   &   & 1 & 8  \\ 
                           &   & 3 & 6 & 4 & 5  \\
                           &   &   & 7 & 9        
\end{array}, \quad
\begin{array}{cccccc} &   &   &   &   & 2  \\  
                           &   &   &   & 1 & 8  \\ 
                           &   & 6 & 7 & 4 & 5  \\
                           &   &   & 9 & 3        
\end{array}
$$
before we finally obtain
$$
\begin{array}{cccccc} &   &   &   &   & 8  \\  
                           &   &   &   & 1 & 5  \\ 
                           &   & 6 & 7 & 4 & 2  \\
                           &   &   & 9 & 3        
\end{array}.$$
Note that the output tabloid is equal to $Q_{\pi^{-1}}$.
This is because in the course of applying $\fj$ to a 
pair $(P_{\pi},P)$ the first tabloid of the current pair can
always be obtained from the second tabloid in the pair by applying $\pi$. 

Observe that the following operation $\bj$ is the inverse of $\fj$. Let
$T'$ be a (shifted) standard skew tableau and $S'$ be a tabloid. For $i=1$ to $i=n$ 
perform backward jeu de taquin in $T'$ in the subshape of $\mu / \lambda$ consisting of the 
cells of $S'$ whose entries are greater or equal than $i$ and
with the entry of $T'$ that is in the cell of the entry $i$ in $S'$. Again perform the transpositions 
simultaneously in $S'$. If $T'$ results in $T$ and $S'$ results in $S$
then we define $\bj(S',T')=(T,S)$. Observe that $S$ is a (shifted) standard skew tableau by
construction. Furthermore $\bj \cdot \fj = \id$ and $\fj \cdot \bj = \id$.

If we apply $\bj$ to  $S'=Q_{\pi^{-1}}$  and $T'=Q$ we obtain the following  pairs 
of intermediate tabloids: 
\begin{eqnarray*}
& \left(\begin{array}{ccccc}    &   &   &   & 2  \\  
                                 &   &   & 1 & 8  \\ 
                                 & 6 & 7 & 4 & 5  \\
                                 &   & 9 & 3        
        \end{array},
   \begin{array}{ccccc}    &   &   &   & 8  \\  
                                 &   &   & 3 & 4  \\ 
                                 & 1 & 2 & 5 & 6  \\
                                 &   & 7 & 9        
        \end{array} \right), \quad
\left(        \begin{array}{ccccc}    &   &   &   & 2  \\  
                                 &   &   & 1 & 8  \\ 
                                 & 3 & 6 & 4 & 5  \\
                                 &   & 7 & 9        
        \end{array},
\begin{array}{ccccc}    &   &   &   & 8  \\  
                                 &   &   & 3 & 4  \\ 
                                 & 9 & 1 & 5 & 6  \\
                                 &   & 2 & 7        
        \end{array} \right),& \\
& \left(\begin{array}{ccccc}    &   &   &   & 2  \\  
                                 &   &   & 1 & 8  \\ 
                                 & 3 & 4 & 6 & 5  \\
                                 &   & 7 & 9        
        \end{array},
   \begin{array}{ccccc}    &   &   &   & 8  \\  
                                 &   &   & 3 & 4  \\ 
                                 & 9 & 5 & 1 & 6  \\
                                 &   & 2 & 7        
        \end{array} \right), \quad
\left(        \begin{array}{ccccc}    &   &   &   & 2  \\  
                                 &   &   & 1 & 5  \\ 
                                 & 3 & 4 & 6 & 8  \\
                                 &   & 7 & 9        
        \end{array},
\begin{array}{ccccc}    &   &   &   & 8  \\  
                                 &   &   & 3 & 6  \\ 
                                 & 9 & 5 & 1 & 4  \\
                                 &   & 2 & 7        
        \end{array} \right)& \\
=(P,R) 
\end{eqnarray*} 

Clearly $\mj_P ( P_\pi) =Q$ is equivalent with $\fj(P_{\pi},P) = (Q_{\pi^{-1}}, Q)$.
Thus the assertion of the theorem is that $\fj$ is an involution. This is equivalent 
to $\fj = \bj$. 

In order to show that we decompose $\fj$ and $\bj$ into its elementary
steps. Let $S$, $T$ be tabloids of the same shape and $1 \le i \le n$.
Then the pair $(T',S')=J^{i}(T,S)$ is defined as follows: Let $T'$ be the
tabloid we obtain by performing  
forward jeu de taquin with the entry in the cell of $T$ which is labelled
with $i$ in $S$ and perform the corresponding transpositions in $S$ also 
in order to obtain $S'$. 
Observe that 
$$\fj(T,S)=\switch J^{1} J^{2} \dots J^{n}(T,S)$$
with $\switch(X,Y)=(Y,X)$ if $T$ is a tabloid and $S$ is a (shifted) standard
skew tableau.
The pair $(S,T)=J_{i}(S',T')$ is defined as 
follows, where $S'$, $T'$ are tabloids of the same shape and the cells of the entries of
$S'$ greater or equal to $i$ form a subshape:
Let $T$ be the tabloid we obtain by performing backward 
jeu de taquin in $T'$ in the subshape of $\mu / \lambda$ consisting of the 
cells of $S'$ whose entries are greater or equal than $i$
and with the entry in the cell of $T'$ which is labelled with $i$ in 
$S'$. Again perform the corresponding transpositions in $S'$ also in order to obtain
$S$.
Observe that 
$$\bj(S',T')=\switch J_{n} J_{n-1} \dots J_{1}(S',T')$$
if $S'$ is a tabloid and $T'$ is a (shifted) standard skew tableau.

\medskip

Next we show the following identity
\begin{equation}
\label{id}
\fj(T,S) = \switch J^{1}  \dots J^{k-1}  J^{k+1}  \dots J^{n}  J_{1}(T,S),
\end{equation}
for a tabloid $T$ and a (shifted) standard skew tableau $S$,
where $k$ is the entry in the cell of $S$ which is labelled with $1$ in $T$.
Let $(i,j)$ be the cell of $1$ in $T$. If the cell $(i,j)$ has 
no neighbour to the left and no neighbour above in $T$ then 
$J_{1}(T,S)=(T,S)$ and since 
$$
J^{k} J^{k+1}  \dots  J^{n}  (T,S)= J^{k+1}  \dots  J^{n} (T,S)
$$
($1$ is always stable),
this proves the assertion in this case.
Now suppose that $(i,j-1),(i-1,j)$ are both cells
in the fixed shape  and $S_{i,j-1} > S_{i-1,j}$ for the other cases are similar. 
Then the entry $1$ in $T$ is first involved in a transposition in the application 
of $\fj$ to $(T,S)$ when performing the first step of $J^{S_{i,j-1}}$ to the current pair. 
In this case the entries in cells $(i,j-1)$ and $(i,j)$ are exchanged in both
tabloids. Note that this is also the first time that an entry in cell $(i,j-1)$ is involved in a
transposition in the application of $\fj$. But
this tranposition is also the first step in the application of 
$J_{1}$ to $(T,S)$. If $(i,j-1)$ has no neighbour to the left and no neighbour
above in the fixed shape the assertion is proved, for $1$ in $T$ is neither
involved in another transposition of $\fj$ nor of $J_{1}$ and the application
of $J_{1}$ terminates. Otherwise suppose 
$(i-1,j-1), (i,j-2)$ are both cells in the fixed
shape and $S_{i-1,j-1} > S_{i,j-2}$ for the other cases are similar. Then the entry $1$ in $T$ is
involved in a transposition in the application of $\fj$ to $(T,S)$ for the
second time when performing the first transposition of $J^{S_{i-1,j-1}}$ to the current
pair. In this step the entries in cells $(i-1,j-1)$ and $(i,j-1)$ are exchanged in
both tabloids. Again this is the first time an entry in $(i-1,j-1)$ is involved in a transposition of
$\fj$.
But this transposition is also the second step in the
application of $J_{1}$ to $(T,S)$ etc. 
Roughly speaking the backward path in the
application of $J_{1}$ to $(T,S)$, which we obtain by performing backward jeu de
taquin with the entry in $S$ in the cell labelled with  $1$ in $T$,
is equal to the 'backward path' of $1$ in $T$, which we
obtain indirectly in the application of $\fj$ to $(T,S)$ by performing forward jeu de 
taquin to all entries in
$T$.   

\medskip

Now we show $\fj=\bj$ by induction with respect to $n$. 
For $n=1$ there is nothing to prove.
Suppose that
$(T',S')$ is a pair of tabloids, where entry $1$ in $T'$ has no neighbour to the left and no
neighbour above in $T'$ and where $S'$ without the entry in the cell of $1$ in 
$T'$ is standard. Let $\fj'(T',S')$ and $\bj'(T',S')$ denote the output pairs 
after the application of 
$\fj$ and $\bj$ to the pair $(T',S')$ with the cell of $1$ in $T'$ omitted in 
both tabloids. By induction $\fj'(T',S')=\bj'(T',S')$.  Let $(T,S)$ be a pair 
of a tabloid $T$ and a (shifted) standard skew tableau $S$ and observe that the 
pair $(T',S') = J_{1} (T,S)$ has the property of $(T',S')$ above. Thus 
$$
\fj (T,S) = \fj' (J_{1} (T,S)) = \bj' (J_{1}(T,S)) = \bj (T,S),
$$
where the first equality follows from (\ref{id}).
\qed

\begin{cor}
Let $P,Q$ be two (shifted) standard skew tableaux of the same shape. Then $A_{P,Q}=A_{Q,P}$.
\end{cor}

{\it Proof.}
By the theorem 
$$
\{ \pi \in {\mathcal S}_n | \mj_{P} (P_\pi)=Q \} = \{ \pi \in {\mathcal S}_n | \mj_{Q} (Q_{\pi^{-1}}) = P \},
$$
where ${\mathcal S}_{n}$ denotes the symmetric group of order $n$.

Thus
\begin{eqnarray*}
A_{P,Q} & = & |\{ T \mbox{ tabloid } | \mj_P(T) = Q \}| \\
        & = & |\{ \pi \in {\mathcal S}_n | \mj_P (P_\pi) = Q \}|  \\
        & = & |\{ \pi \in {\mathcal S}_n | \mj_Q (Q_{\pi^{-1}}) = P \}| \\ 
        & = & |\{ T  \mbox{ tabloid } | \mj_Q(T) = P \}| \\
        & = & A_{Q,P}.
\end{eqnarray*}
\qed

\end{document}